\documentclass{kms-b}

\issueinfo{}
  {}
  {}
  {}
\pagespan{1}{}
\copyrightinfo{}

\usepackage{graphicx}
\allowdisplaybreaks

\theoremstyle{plain}

\theoremstyle{definition}

\theoremstyle{remark}

\begin{document}

\title[Estimating the distance Estrada index]
{Estimating the distance Estrada index}

\author[Y. Shang]{Yilun Shang}
\address{Yilun Shang\\ Einstein Institute of Mathematics \\ Hebrew University \\ Jerusalem 91904, Israel}
\email{shylmath@hotmail.com}



\begin{abstract}
Suppose $G$ is a simple graph on $n$ vertices. The $D$-eigenvalues
$\mu_1,\mu_2,\cdots,\mu_n$ of $G$ are the eigenvalues of its
distance matrix. The distance Estrada index of $G$ is defined as
$DEE(G)=\sum_{i=1}^ne^{\mu_i}$. In this paper, we establish new
lower and upper bounds for $DEE(G)$ in terms of the Wiener index
$W(G)$. We also compute the distance Estrada index for some concrete
graphs including the buckminsterfullerene $C_{60}$.

\end{abstract}

\maketitle

\section{Introduction}

Let $G$ be a simple $n$-vertex graph with vertex set $V(G)$. Denote
by $D(G)=(d_{ij})\in\mathbb{R}^{n\times n}$ the distance matrix of
$G$, where $d_{ij}$ signifies the length of shortest path between
vertices $v_i\in V(G)$ and $v_j\in V(G)$. Then the adjacency matrix
$A(G)=(a_{ij})$ of the graph can be defined by $a_{ij}=1$ if
$d_{ij}=1$, and $a_{ij}=0$ otherwise. Since $D(G)$ is a real
symmetric matrix, its eigenvalues are real numbers. We order the
eigenvalues in a non-increasing manner as
$\mu_1\ge\mu_2\ge\cdots\ge\mu_n$ (they are customarily called
$D$-eigenvalues of $G$ \cite{2}). The distance Estrada index of $G$
is defined as
\begin{equation}
DEE(G)=\sum_{i=1}^ne^{\mu_i}.\label{1}
\end{equation}

This graph-spectrum-based structural invariant is recently proposed
in \cite{1}, and some results on its bounds can be found in
\cite{1a,1b,a}. If we replace in (\ref{1}) the $D$-eigenvalues
$\{\mu_i\}_{i=1}^n$ by the eigenvalues $\{\lambda_i\}_{i=1}^n$ of
the adjacency matrix $A(G)$, we recover the well-researched graph
descriptor Estrada index \cite{3}. The Estrada index can be used as
an efficient measuring tool in a number of areas in chemistry and
physics, and its mathematical properties have been intensively
studied (see e.g. \cite{10,7,4,5,6,9,11,11a,m1,m2,8}, to mention
only a few).

Apart from its formal analogy to the Estrada index, we believe the
distance Estrada index (\ref{1}) is potentially of vast importance
in physical chemistry. After all, the most natural description of a
molecular graph is in terms of the distances---be them geometric or
topological---between pairs of vertices. The oldest distance-based
invariant, perhaps, is the Wiener index \cite{12}, which has found
useful applications in structure---property correlations; see e.g.
\cite{13,14,16,15}. In this paper, we aim to establish lower and
upper bounds for $DEE(G)$ by using the Wiener index, which allows us
to gain insight into the relationship between the distance Estrada
index and the Wiener index, and, in particular, gain better
understanding of the dependence of the distance Estrada index on the
concept of distance degree whereby the Wiener index is constructed.
We mention that various properties of $D(G)$ for some interesting
graphs can be found in e.g. \cite{a3,a2}.

The rest of the paper is organized as follows. In Section 2, we give
some notations and lemmas. The bounds for $DEE(G)$ are provided in
Section 3. In Section 4, we compute the distance Estrada index of
some concrete graphs, including the buckminsterfullerene $C_{60}$,
to demonstrate the availability of our obtained results.

\section{Preliminaries}

Let $G$ be a simple connected $n$-vertex graph with vertex set
$V(G)=\{v_1,v_2,$ $\cdots,v_n\}$. Denote by $D(G)$ the distance
matrix of the graph $G$. The distance degree of a vertex $v_i$ is
given by $D_i=\sum_{j=1}^nd_{ij}$. This concept first appears in
\cite{17} and is reinvented recently in \cite{18} under the name of
distance degree. The Wiener index \cite{12} of $G$, denoted by
$W(G)$, is the sum of the distances between all (unordered) pairs of
vertices of $G$, that is
\begin{equation}
W(G)=\sum_{i<j}d_{ij}=\frac12\sum_{i=1}^nD_i.\label{2}
\end{equation}
Let $M(G)=(\prod_{i=1}^nD_i)^{1/n}$ be the geometric mean of the
distance degrees. Then $2W(G)/n\ge M(G)$ holds, and equality is
attained if and only if $D_1=D_2=\cdots=D_n$ (i.e., the graph $G$ is
distance degree regular \cite{17}).

Several properties of the spectrum of the distance matrix $D(G)$
follows easily from its definition. For $k\in\mathbb{N}$, let
$N_k=\sum_{i=1}^n\mu_i^k$ be the $k$th spectral moment. Since all
elements of $D(G)$ are integers, all moments $N_k$ are also
integral. In particular, $N_1=0$, i.e., $D(G)$ is traceless; and
$N_2=2\sum_{i<j}d_{ij}^2$. The following two lemmas will be needed
later.

\smallskip
\noindent\textbf{Lemma 1.} \cite{18}\itshape \quad A connected graph
$G$ has two distinct $D$-eigenvalues if and only if $G$ is a
complete graph. \normalfont
\smallskip

\smallskip
\noindent\textbf{Lemma 2.} \cite{19}\itshape \quad Let
$a_1,a_2,\cdots,a_n$ be nonnegative numbers. Then
$$
n\left(\frac1n\sum_{i=1}^na_i-\left(\prod_{i=1}^na_i\right)^{\frac1n}\right)\le
n\sum_{i=1}^na_i-\left(\sum_{i=1}^na_i^{\frac12}\right)^2.
$$
\normalfont
\smallskip

\section{Bounds for the distance Estrada index}

Our main result reads as follows.

\smallskip
\noindent\textbf{Theorem 1.} \itshape \quad Let $G$ be a connected
graph on $n$ vertices. Denote by $\Delta(G)$ the diameter of $G$.
Then
\begin{equation}
e^{\left(\frac{4W^2(G)-M^2(G)n}{n(n-1)}\right)^{\frac12}}+\frac{n-1}{e^{\frac{1}{n-1}\left(\frac{4W^2(G)-M^2(G)n}{n(n-1)}\right)^{\frac12}}}
\le DEE(G)\le
n-1+e^{\sqrt{2}\Delta^{\frac12}(G)W^{\frac12}(G)}.\label{3}
\end{equation}
The equality on the left-hand side of (\ref{3}) holds if and only if
$G$ is the complete graph $K_n$. The equality on the right-hand side
of (\ref{3}) holds if and only if $G=K_1$, i.e., a single vertex.
\normalfont
\smallskip

\noindent\textbf{Remark 1.} When $n=1$, we will have
$4W^2(G)=M^2(G)n$ and view the leftmost term of (\ref{3}) as 1.

\noindent\textbf{Proof}. \textit{Lower bound.} Using the
arithmetic-geometric mean inequality, we obtain
\begin{eqnarray}
DEE(G)&=&\sum_{i=1}^ne^{\mu_i}\nonumber\\
&\ge&e^{\mu_1}+(n-1)\left(\prod_{i=2}^ne^{\mu_i}\right)^{\frac{1}{n-1}}\nonumber\\
&=&e^{\mu_1}+(n-1)e^{-\frac{\mu_1}{n-1}},\label{4}
\end{eqnarray}
where we have used the fact that $N_1=\sum_{i=1}^n\mu_i=0$.

In \cite{18} it was shown that
\begin{equation}
\mu_1\ge\left(\frac{\sum_{i=1}^nD_i^2}{n}\right)^{\frac12}.\label{5}
\end{equation}
Setting $\sqrt{a_i}=D_i$ in Lemma 2, we get
$$
n^2\left(\frac{\sum_{i=1}^nD_i^2}{n}-\left(\frac{2W}{n}\right)^2\right)\ge\sum_{i=1}^nD_i^2-n\left(\prod_{i=1}^nD_i^2\right)^{\frac1n}.
$$
Combining this with (\ref{5}) yields
\begin{equation}
\mu_1\ge\left(\frac{4W^2(G)-M^2(G)n}{n(n-1)}\right)^{\frac12}\ge0.\label{6}
\end{equation}
Clearly, $4W^2(G)=M^2(G)n$ (namely, the second equality in (\ref{6})
holds) if and only if $n=1$.

It is elementary to show that for $n\ge1$ the function
\begin{equation}
f(x)=e^x+\frac{n-1}{e^{\frac{x}{n-1}}}\label{func}
\end{equation}
monotonically increases in the interval $[0,+\infty)$ (here we take
the limit function $f(x)=e^x$ when $n=1$). Therefore, by means of
(\ref{4}) and (\ref{6}) we arrive at the first half of Theorem 1.

Note that when $G=K_n$, we have $\mu_1=n-1$,
$\mu_2=\cdots=\mu_n=-1$, $W(G)=n(n-1)/2$ and $M(G)=n-1$. Hence,
$DEE(G)=e^{n-1}+(n-1)e^{-1}$ and the equality on the left-hand side
of (\ref{3}) holds. Conversely, suppose that the equality holds,
then from (\ref{4}) we have $\mu_2=\cdots=\mu_n$. We assume that
$n\ge2$. It follows from (\ref{6}) that $\mu_1>0$. Then $G$ has
exactly two distinct $D$-eigenvalues, and Lemma 1 indicates that $G$
is the complete graph $K_n$.

\textit{Upper bound.} Let $n_+$ be the number of positive
$D$-eigenvalues of $G$, we obtain
\begin{eqnarray*}
DEE(G)&\le&n-n{_+}+\sum_{i=1}^{n_+}e^{{\mu_i}}\\
&=&n-n{_+}+\sum_{i=1}^{n_+}\sum_{k=0}^{\infty}\frac{\mu_i^k}{k!}\\
&=&n+\sum_{k=1}^{\infty}\frac{1}{k!}\sum_{i=1}^{n_+}\mu_i^k\\
&\le&n+\sum_{k=1}^{\infty}\frac{1}{k!}\left(\sum_{i=1}^{n_+}\mu_i^2\right)^{\frac{k}{2}}\\
&=&n+\sum_{k=1}^{\infty}\frac{1}{k!}\left(2\sum_{i<j}d_{ij}^2-\sum_{i=n_++1}^n\mu_i^2\right)^{\frac{k}{2}}\\
&\le&n-1+e^{\sqrt{2\sum_{i<j}d_{ij}^2}}\\
&\le&n-1+e^{\sqrt{2\Delta\sum_{i<j}d_{ij}}},
\end{eqnarray*}
which directly leads to the right-hand side inequality in (\ref{3}).

From the above derivation it is apparent that equality holds if and
only if the graph $G$ has all zero $D$-eigenvalues. Since $G$ is a
connected graph, this only happens when $G=K_1$ (and thus
$DEE(G)=1$).

The proof of Theorem 1 is completed. $\Box$

\smallskip
\noindent\textbf{Remark 2.} In \cite{1b}, it was proved that
\begin{equation}
DEE(G)\ge
e^{\frac{2W(G)}{n}}+(n-1)e^{-\frac{2W(G)}{n(n-1)}}.\label{7}
\end{equation}
If we utilize the property $2W(G)/n\ge M(G)$, then we obtain
$$
\frac{2W(G)}{n}\le\left(\frac{4W^2(G)-M^2(G)n}{n(n-1)}\right)^{\frac12}.
$$
Since the function $f(x)$ defined in (\ref{func}) is strictly
increasing, we see that our lower bound in (\ref{3}) is better than
the bound in (\ref{7}).

\smallskip
\noindent\textbf{Remark 3.} It was shown in \cite{1} that
\begin{equation}
DEE(G)\le n-1+e^{\Delta(G)\sqrt{n(n-1)}}.\label{8}
\end{equation}
Since $d_{ij}\le\Delta(G)$ for all $i$ and $j$,
$n(n-1)\Delta(G)\ge2W(G)$. Obviously, our upper bound in (\ref{3})
is better than the bound in (\ref{8}).

If the graph $G$ is $r$-distance regular for some $r\in\mathbb{N}$,
we have $D_1=D_2=\cdots=D_n=r$ \cite{18}. Consequently, $W(G)=nr/2$
and $M(G)=r$. The following result is immediate.

\smallskip
\noindent\textbf{Corollary 1.} \itshape \quad Let $G$ be a connected
$r$-distance regular graph on $n$ vertices. Denote by $\Delta(G)$
the diameter of $G$. Then
\begin{equation}
e^r+\frac{n-1}{e^{\frac{r}{n-1}}} \le DEE(G)\le
n-1+e^{\sqrt{\Delta(G)nr}}.\label{3a}
\end{equation}
The equality on the left-hand side of (\ref{3a}) holds if and only
if $G$ is the complete graph $K_n$ with $n=r+1$. The equality on the
right-hand side of (\ref{3a}) holds if and only if $G=K_1$, i.e., a
single vertex. \normalfont
\smallskip

\section{Some examples}

In this section, we provide some concrete examples to demonstrate
the calculations of distance Estrada index as well as the
feasibility of the above obtained results.

\smallskip
\noindent\textbf{Example 1.} In this example, the graph $G$ is the
cycle over $n=6$ vertices, namely, a hexagonal cell. Its distance
matrix is shown below
$$
D(G)=\left(\begin{array}{cccccc}
0&1&2&3&2&1\\
1&0&1&2&3&2\\
2&1&0&1&2&3\\
3&2&1&0&1&2\\
2&3&2&1&0&1\\
1&2&3&2&1&0
\end{array}\right).
$$
Since $D(G)$ is a circulant matrix, the $D$-eigenvalues consist of
$\{\omega_j+2\omega_j^2+3\omega_j^3+2\omega_j^4+\omega_j^5\}_{j=0}^5$,
where $\omega_j=\cos(2\pi j/n)+i\sin(2\pi j/n)$, $i=\sqrt{-1}$, and
$j=0,1,\cdots,5$. Via some simplifications, we have
$\mu_1=9,\mu_2=\mu_3=0,\mu_4=-1$ and $\mu_5=\mu_6=-4$. Thus, we
calculate that $DEE(G)=\sum_{k=1}^6e^{\mu_k}=8105.5$.

It is evident that $G$ is a connected $r$-distance regular graph
with $r=9$. The diameter is given by $\Delta(G)=3$. Therefore, from
Corollary 1 we have the following bounds
$$
8103.9\le DEE(G)\le 337033.2.
$$
It turns out that the lower bound is very sharp while the upper
bound is conservative in this specific example.

\begin{figure}[!t]
\begin{center}
\scalebox{0.5}{\includegraphics[138pt,203pt][456pt,354pt]{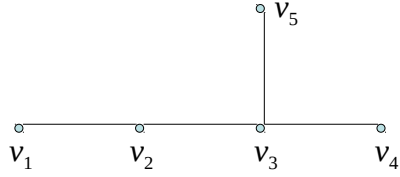}}\caption{
A chemical tree $G$ on vertex set $V(G)=\{v_1,v_2,\cdots,v_5\}$.
Recall that a chemical tree is a tree having no vertex with degree
greater than 4.}
\end{center}
\end{figure}

\smallskip
\noindent\textbf{Example 2.} In Fig. 1 we display a chemical tree of
$n=5$ vertices. The distance matrix is
$$
D(G)=\left(\begin{array}{ccccc}
0&1&2&3&3\\
1&0&1&2&2\\
2&1&0&1&1\\
3&2&1&0&2\\
3&2&1&2&0
\end{array}\right).
$$
The $D$-eigenvalues of the graph are as follows:
$\mu_1=7.46,\mu_2=-0.51,\mu_3=-1.08,\mu_4=-2$ and $\mu_5=-3.86$. We
obtain $DEE(G)=1738.2$.

On the other hand, we can easily obtain $W(G)=18$, $M(G)=7.04$ and
$\Delta(G)=3$. Hence, from Theorem 1 we have the following bounds
$$
1393.4\le DEE(G)\le 32611.7.
$$

\begin{figure}[!t]
\begin{center}
\scalebox{0.6}{\includegraphics[0pt,0pt][184pt,173pt]{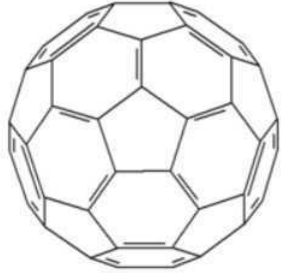}}\caption{A
projection depiction of the buckminsterfullerene $C_{60}$.}
\end{center}
\end{figure}

\smallskip
\noindent\textbf{Example 3.} In this example, we consider the
buckminsterfullerene $C_{60}$, which is a well-known member of the
fullerene family \cite{20}. As a graph, $C_{60}$ is a truncated
icosahedron with $n=60$ vertices and 32 faces (including 20 hexagons
and 12 pentagons); see Fig. 2 for an illustration. The
$D$-eigenvalues of $C_{60}$ were computed in \cite{21} by using the
Givens-Householder method (see Table 1 in \cite{21}). For example,
there are exactly 18 positive $D$-eigenvalues and no zero
$D$-eigenvalue. Based on these $D$-eigenvalues we calculate by
definition (\ref{1}) that $DEE(C_{60})=152.11+e^{278}$.

Since $C_{60}$ contains 60 vertices, a wise approach to capture its
distance matrix is to understand the associated distance level
diagram, as is done in Figure 1 of \cite{21}. From that, we conclude
that $C_{60}$ is a connected $r$-distance regular graph with $r=278$
and the diameter is $\Delta(C_{60})=9$. Hence, $W(C_{60})=nr/2=8340$
and $M(C_{60})=r=278$. Corollary 1 leads to the following bounds
$$
0.53+e^{278}\le DEE(C_{60})\le 59+e^{387}.
$$
Clearly, the upper bound obtained is way above the true value for
$DEE(C_{60})$. Future work may study better behaved estimations.

\end{document}